# THE UBIQUITY OF THOMPSON'S GROUP $F$ IN GROUPS OF PIECEWISE LINEAR HOMEOMORPHISMS OF THE UNIT INTERVAL[1]

MATTHEW G. BRIN

Department of Mathematical Sciences

State University of New York at Binghamton

Binghamton, NY 13902-6000

USA

April 13, 1997

## 0. INTRODUCTION

In the 1960s, Richard J. Thompson introduced a triple of groups $F \subseteq T \subseteq G$ which, among them, supplied: the first examples of infinite, finitely presented, simple groups [14] (see [6] for published details); a technique for constructing an elementary example of a finitely presented group with an unsolvable word problem [12]; the universal obstruction to a problem in homotopy theory [8]; the first examples of torsion free groups of type $FP_\infty$ and not of type $FP$ [5]. In abstract measure theory, it has been suggested by Geoghegan (see [3] or [9, Question 13]) that $F$ might be a counterexample to the conjecture that any finitely presented group with no non-cyclic free subgroup is amenable (admits a bounded, non-trivial, finitely additive measure on all subsets that is invariant under left multiplication). Recently, $F$ has arisen in the theory of groups of diagrams over semigroup presentations [10], and in the algebra of string rewriting systems [7]. For more extensive bibliographies and more results on Thompson's groups and their generalizations see [1, 4, 6].

A persistent peculiarity of Thompson's groups is their ability to pop up in diverse areas of mathematics. This suggests that there might be something very natural about Thompson's groups. We support this idea by showing (Theorem 1 below) that $PL_o(I)$, the group of piecewise linear (finitely many changes of slope), orientation preserving, self homeomorphisms of the unit interval, is riddled with copies of $F$: a very weak criterion implies that a subgroup of $PL_o(I)$ must contain an isomorphic copy of $F$.

## 1. STATEMENT

Throughout the paper $I$ will represent the unit interval $[0, 1]$, functions will be written to the right of their arguments, and composition will proceed from left to right.

---

[1] AMS Classification (1991): primary 20F05

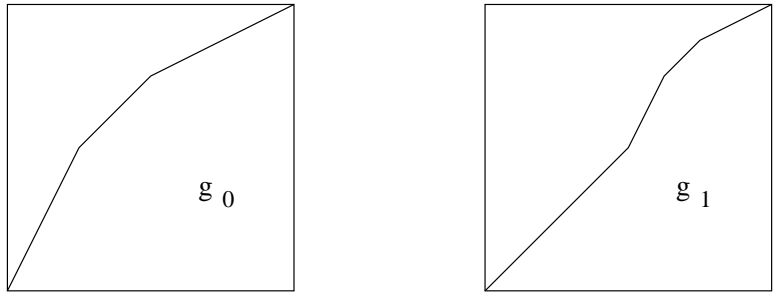

FIGURE 1. $f_0$ and $f_1$

We start with two equivalent presentations that define $F$ and give a reference for the equivalence. We then give a standard embedding of $F$ in $PL_o(I)$ that is sometimes used as its definition. The embedding gives a good picture of $F$. We use the embedding to point out the well known fact that $F$ itself is riddled with copies of $F$. We then state our result. The proof follows in the next section.

We define $F$ by the presentations

(1) $\qquad F = \langle g_0, g_1 \mid g_1^{g_0 g_0} = g_1^{g_0 g_1}, g_1^{g_0 g_0 g_0} = g_1^{g_0 g_0 g_1} \rangle, \qquad$ and

(2) $\qquad F = \langle g_0, g_1, g_2, \ldots \mid g_j^{g_i} = g_{j+1}, \text{ for all } i < j \rangle$

where $a^b = b^{-1}ab$. Note that $g_i = g_0^{1-i} g_1 g_0^{i-1}$ for all $i > 1$ in (2). There is an isomorphism that carries the pair $(g_0, g_1)$ in (1) to the pair in (2) with the same names [6, Thm. 3.1]. The relations in (1) above and for $F_1$ in [6] are identical in spite of differences in appearance.

There is an isomorphism from $F$ to the group $\bar{F}$ of elements of $PL_o(I)$ that have slopes that are integral powers of 2, and that have "breaks" (discontinuities of slope) restricted to $\mathbf{Z}[\frac{1}{2}]$ (the set of all $a/2^b$ with $a$ and $b$ in the integers $\mathbf{Z}$). The isomorphism carries the pair $(g_0, g_1)$ to $(f_0, f_1)$ defined by

$$xf_0 = \begin{cases} 2x & 0 \le x \le \frac{1}{4}, \\ x + \frac{1}{4} & \frac{1}{4} \le x \le \frac{1}{2}, \\ \frac{1}{2}x + \frac{1}{2} & \frac{1}{2} \le x \le 1, \end{cases} \qquad xf_1 = \begin{cases} x & 0 \le x \le \frac{1}{2}, \\ 2x & \frac{1}{2} \le x \le \frac{5}{8}, \\ x + \frac{1}{8} & \frac{5}{8} \le x \le \frac{3}{4}, \\ \frac{1}{2}x + \frac{1}{2} & \frac{3}{4} \le x \le 1. \end{cases}$$

The graphs of $f_0$ and $f_1$ are shown in Figure 1.

Conjugating by $x \mapsto x/2$ gives an isomorphism from $\bar{F}$ to its subgroup of elements that are fixed off $[0, \frac{1}{2}]$, and conjugating further by $x \mapsto x + \frac{1}{4}$ gives an isomorphism from $\bar{F}$ to its subgroup of elements fixed off $[\frac{1}{4}, \frac{3}{4}]$. For $0 < a < b < 1$ with $a, b \in \mathbf{Z}[\frac{1}{2}]$, there is an element $f \in \bar{F}$ carrying $[\frac{1}{4}, \frac{3}{4}]$ to $[a, b]$ (see e.g., [6, Lemma 4.2] or [1, Lemma 2.1]), and $\bar{F}$ is thus isomorphic to its subgroup of elements fixed off $[a, b]$. We see that $F$ contains many copies of itself.

It is an easy exercise that carrying $(g_0, g_1)$ to $(f_0, f_1)$ preserves the relations of (1), giving a homomorphism from $F$ to $\bar{F}$. The group presented by (2) has no proper non-abelian quotients [6, Thm. 4.3], [8, (L12)] or [2, Lemma 2.1.5] and $f_0$ and $f_1$ do not commute, so the homomorphism is a monomorphism. Then [6, Cor. 2.6] or [4, prop. 4.4] shows that the monomorphism is an epimorphism. (Although we give recent references, most of this was recorded in [14]. Note that the functions $A$ and $B$ in [6] are the inverses of $f_0$ and $f_1$ and composition is from right to left. Thus relations in [6] read similarly to the relations here.)

We give definitions needed in the statement of Theorem 1. In the definitions, $H$ will be a subgroup of $PL_o(I)$.

The *fixed point set* $F(H)$ of $H$ is the set points fixed by all $h \in H$. The *support* of $H$ is the complement of $F(H)$ in $I$. The *orbit* of an $x \in I$ under $H$ is the set $\{xh \mid h \in H\}$. The orbit of an $x \in F(H)$ is a single point, and it is elementary that the orbit of an $x \notin F(H)$ is infinite. For an $x \notin F(H)$, the *orbital* of $x$ under $H$ is the convex hull of the orbit of $x$ under $H$ and is the set of all $y \in I$ for which there are $h_1$ and $h_2$ in $H$ with $xh_1 < y < xh_2$. It is easy to see that the orbitals of $H$ are the components of the support of $H$. If $h \in H$, then $\langle h \rangle$ is the cyclic subgroup of $H$ generated by $h$ and the orbitals of $\langle h \rangle$ will be referred to as the *orbitals* of $h$. The *support* of $h$ is the support of $\langle h \rangle$.

If $h$ is in $H$ and $A$ is an orbital of $H$ with an endpoint $a$, then we say that $h$ *approaches* $a$ in $A$ if some orbital $B$ of $h$ lies in $A$ and has $a$ as an endpoint.

**Theorem 1.** *Let $H$ be a subgroup of $PL_o(I)$. Assume that $H$ has an orbital $(a, b)$ and that some element of $H$ approaches one end of $a$ or $b$ in $(a, b)$ but not the other. Then $H$ contains a subgroup isomorphic to $F$.*

Thompson's groups were generalized by Higman in [11] and Brown in [4] to three infinite families of groups ($F_n$, $T_{n,r}$ and $G_{n,r}$ with $F = F_2$, $T = T_{2,1}$ and $G = G_{2,1}$). Arguments similar to ours would show that $PL_o(I)$ has many copies of the $F_n$. However [2], the $F_n$ are isomorphic to subgroups of $F$, so separate arguments for the $F_n$ are not needed.

## 2. Arguments

Theorem 1 is proven by looking at one orbital of $H$ at a time. Our main technical lemma finds copies of $F$ in single orbitals. We need some definitions.

If $H$ is a subgroup of $PL_o(I)$ and $A$ is an orbital of $H$, then $H_A$ will denote the restrictions to $A$ of all elements of $H$. The elements of $H_A$ form a group which is a homomorphic image of $H$. We call $H_A$ the *projection* of $H$ onto the orbital $A$.

If $H$, $h$, $A$, $a$ and $B$ are as in the defintion that $h$ approaches the endpoint $a$ of $A$ in $A$, then either $xh^i$ limits to $a$ as $i$ goes to $+\infty$ for all $x \in B$, or $xh^i$ limits to $a$ as $i$ goes to $-\infty$ for all $x \in B$. In the first case, we say that $h$ approaches $a$ *positively* in $A$ and define the *sign of the approach* of $h$ to $a$ in $A$, written $\sigma(h, a, A)$, to be $+1$. In the second case, we say that $h$ approaches $a$ *negatively* in $A$ and set $\sigma(h, a, A) = -1$.

The main technical lemma follows. The statement is given in two parts to split the principle conclusion from the technicalities needed in the proof of Theorem 1.

**Lemma 2.1.** *Part I: Let $H$ be a subgroup of $PL_o(I)$. Assume that $A = (a, b)$ is the only orbital of $H$. Let $e$ be one of $a$ or $b$ and let $e'$ be the other. Assume that some element $h \in H$ approaches $e$ and not $e'$. Then there is a subgroup of $H$ isomorphic to $F$.*

*Part II: The isomorphism between the subgroup of $H$ and $F$ is obtained by carrying a pair of elements $(g_0, g_1)$ of $H$ to the similarly named pair in Presentation (1), where $(g_0, g_1)$ is as specified in (i) or (ii) below.*

(i) *There are $f_1 \in H$ and $K > 0$ so that if $\sigma = \sigma(h, e, A)$, then we can set $(g_0, g_1) = (h^{i\sigma}, (f_1^{-1} h f_1)^{i\sigma})$ for any $i > K$.*

(ii) *If there is an $f \in H$ with $A$ as its orbital and with $f = h$ near $e$ in $A$, then there exists $K > 0$ so that if $\sigma = \sigma(f, e, A) = \sigma(h, e, A)$, then we can set $(g_0, g_1) = (f^{i\sigma}, h^{i\sigma})$ for any $i > K$.*

*In both (i) and (ii), the support of $g_1$ is a proper subset of one orbital of $g_0$.*

*Proof.* We will give the proof assuming $e = b$ and leave it to the reader to check that the same proof applies if $e = a$ with only minor adjustments of some inequalities.

To reduce item (i) to item (ii), we let $B$ be the orbital of $h$ that ends at $b$. There is an $f_1$ in $H$ that carries the greatest lower bound of the support of $h$ into $B$. Thus $f_1^{-1} h f_1$ will have support in $B$, will approach $b$ but not the other end of $B$, and will have slope near $b$ in $B$ the same as the slope of $h$ near $b$ in $B$. We now replace $H$ by the group $G$ generated by $h$ and $f_1^{-1} h f_1$ and note that $G$ is isomorphic to its projection $G_B$ to $B$ since both of the generators are non-trivial on $B$ and one of the generators is trivial off $B$. Thus we make a second replacement of $G$ by its projection $G_B$ to $B$ and the hypotheses of item (ii) are now satisfied with $B$ playing the role of $A$, with $h$ playing the role of $f$ and with $f_1^{-1} h f_1$ playing the role of $h$. Thus item (i) and the last sentence in the statement will follow if we prove item (ii). We now assume that the hypotheses of item (ii) apply to $H$.

Since $h$ and $f$ are PL on a closed interval, they are affine on a neighborhood of $b$ in $A$. Our assumptions dictate that they will have the same slope $\lambda \neq 1$ on that interval. Let $w$ be the greatest lower bound of the support of $h$, and let $x$ be the least so that both $h$ and $f$ are affine on $(x, b)$. We have $a < w < x < b$.

There is a positive integer $m$ so that $y = wf^{m\sigma} > x$. This implies that the support of $q = f^{-m\sigma}hf^{m\sigma}$ lies in $(y, b) \subseteq (x, b)$, an interval on which both $f$ and $h$ are affine with slope $\lambda$. Since $q$ is a conjugate of $h$ it is also affine with slope $\lambda$ on some $(z, b)$ with $z < b$. There is a positive integer $n$ so that $yh^{n\sigma} > z$. We let $K$ be the larger of $m$ and $n$.

Now fix an $i > K$ and define $g_0 = f^{i\sigma}$ and $g_1 = h^{i\sigma}$. Both have fixed points at $b$ and have slope $\lambda^{i\sigma}$ on $(x, b)$. Let $\Lambda$ be the unique affine function with slope $\lambda$ and fixed point $b$. Let

$$\begin{aligned}
g_2 &= g_0^{-1} g_1 g_0 \\
&= f^{-i\sigma} h^{i\sigma} f^{i\sigma} \\
&= f^{-(i-m)\sigma}(f^{-m\sigma} h^{i\sigma} f^{m\sigma}) f^{(i-m)\sigma} \\
&= f^{-(i-m)\sigma} q^{i\sigma} f^{(i-m)\sigma} \\
&= \left( f^{-(i-m)\sigma} q f^{(i-m)\sigma} \right)^{i\sigma} \\
&= \left( \Lambda^{-(i-m)\sigma} q \Lambda^{(i-m)\sigma} \right)^{i\sigma}
\end{aligned}$$

since $q$ has support in $(y, b) \subseteq (x, b)$ on which $f = \Lambda$. Further, the support of $g_2$ lies in $(y, b)\Lambda^{(i-m)\sigma} \subseteq (x, b)$ since positive powers of $\Lambda^\sigma$ move points toward $b$. We also have that $g_2$ is affine on $(z, b)\Lambda^{(i-m)\sigma}$ on which it agrees with $\Lambda^{i\sigma}$.

If we now set $g_3$ to be the conjugate of $g_2$ by $\Lambda^{i\sigma}$, then we also have that $g_3$ is the conjugate of $g_2$ by $g_1 = h^{i\sigma}$ and the conjugate of $g_2$ by $g_0 = f^{i\sigma}$ since $\Lambda$, $h$ and $f$ all agree on the support of $g_2$. The support of $g_3$ lies in

$$(y, b)\Lambda^{(i-m)\sigma}\Lambda^{i\sigma} = (y, b)\Lambda^{i\sigma}\Lambda^{(i-m)\sigma} \subseteq (z, b)\Lambda^{(i-m)\sigma}$$

since $yh^{i\sigma} > z$ and $h = \Lambda$ on $(x, b)$ which contains $y$. Now the conjugates of $g_3$ by each of $g_2$, $g_1$ and $g_0$ are all equal since all of these conjugators equal $\Lambda^{i\sigma}$ on the support of $g_3$. This finishes verifying the relations of Presentation (1). The induced homomorphism from $F$ as presented in (1) to the subgroup of $H$ generated by $g_0$ and $g_1$ is an isomorphism since $g_0$ and $g_1$ do not commute (conjugating $g_1$ by $g_0$ changes the support of $g_1$) and $F$ has no proper non-abelian quotients ([6, Thm. 4.3], [8, (L12)] or [2, Lemma 2.1.5]). □

We now discuss the course of the proof of Theorem 1. Recall the statement.

**Theorem 1.** *Let $H$ be a subgroup of $PL_o(I)$. Assume that $H$ has an orbital $(a, b)$ and that some element of $H$ approaches one end of $a$ or $b$ in $(a, b)$ but not the other. Then $H$ contains a subgroup isomorphic to $F$.*

The hypothesis and Lemma 2.1 immediately give that a quotient of $H$ contains a subgroup isomorphic to $F$. However, this does not immediately imply that $H$ does. The proof will proceed by repeatedly taking subgroups of $H$ that also satisfy the hypotheses of the theorem and that satisfy an accumulated sequence of extra properties. The extra properties will relate to the behavior of the subgroups on orbitals of $H$ other than the orbital $(a, b)$. In what follows, we will give some preliminary facts and lemmas and explain how the proof will proceed. Then the proof will be given as a sequence of steps that close in on the desired subgroup.

We start with some remarks, definitions and lemmas that will be useful in the proof. In the remainder of this section, $H$ will always denote a subgroup of $PL_o(I)$. Unless explicitly stated, the hypotheses of Theorem 1 will not be assumed.

**Remark 2.2.** *If $G \subseteq H$, then each orbital of $G$ is contained in some orbital of $H$.*

The number of orbitals of $G$ might exceed the number of orbitals of $H$.

**Remark 2.3.** *If $G \subseteq H$, $A$ is an orbital of $H$ for which $H_A$ is abelian, and $B$ is an orbital of $G$ for which $B \subseteq A$, then $G_B$ is abelian.*

Let $A_1$ be the orbital $(a, b)$ in the hypothesis of Theorem 1. From Lemma 2.1, we know that two elements $g_0$ and $g_1$ of $H$ have images in $H_{A_1}$ that generate a copy of $F$. If we replace $H$ by the group generated by $g_0$ and $g_1$, then the details of Lemma 2.1 imply that the hypothesis of Theorem 1 is still satisfied. This improves $H$ in that it is now known to be finitely generated and has only finitely many orbitals. For one of its orbitals (we still call it $A_1$), we know that $H_{A_1}$ is isomorphic to $F$. We will be troubled by the projections of $H$ to its other orbitals. We also know that $H$ is a subgroup of the direct product of the projections of $H$ to its various orbitals. Because of the next lemma, our strategy will be to modify $H$ until its projections to its orbitals other than $A_1$ are abelian or other copies of $F$.

**Lemma 2.4.** *Let a non-abelian group $H$ be a subgroup of $H_1 \times H_2$ generated by elements $g_0$ and $g_1$. Assume that $H_2$ is abelian and that the images of the $g_0$ and $g_1$ in $H_1$ satisfy the relations of Presentation (1). Then there is an isomorphism from $F$ to $H$ carrying each $g_i$ of Presentation (1) to $g_i$ in $H$.*

*Proof.* This follows because the relations in Presentation (1) are satisfied by elements that commute, and because $F$ has no proper non-abelian quotients. □

We need to classify orbitals with respect to certain generators. Let $H$ be generated by a fixed ordered pair $(h_0, h_1)$. Let $A$ be an orbital of $H$. We say that $A$ is of *type* $(m_0, m_1)$ where $m_i$ is the number of endpoints of $A$ approached by $h_i$ in $A$. Since $A$ has two endpoints, it is not possible for $m_0 + m_1$ to be less than 2. Thus the possible types are $(1, 1)$, $(2, 0)$, $(0, 2)$, $(2, 1)$, $(1, 2)$ and $(2, 2)$. A consequence of Lemma 2.1 is that $H_A$ contains a subgroup isomorphic to $F$ if $A$ is of type $(1,1)$, $(1,2)$ or $(2,1)$.

With $H$ and $(h_0, h_1)$ as above, we say that an orbital $A$ of $H$ is of *pure type* $(2,1)$ if $A$ is an orbital of $h_0$. Similarly, $A$ is of pure type $(1,2)$ if $A$ is an orbital of $h_1$. We say that an orbital $A$ of $H$ is of pure type $(2,0)$ if $A$ is an orbital of $h_0$, and in addition we say that it is *empty* if the restriction of $h_1$ to $A$ is the identity. Symmetric definitions are made for orbitals of type $(0,2)$. If $A$ of type $(2,0)$ or $(0,2)$ is pure and empty, then the projection of $H$ to $A$ is cyclic and thus abelian.

Orbitals of pure type $(2,1)$ will be special (Lemma 2.1 guarantees at least one), and also a problem in that there are several subtypes that will not cooperate. We need notation to identify the subtypes. Let $A$ be an orbital of pure type $(2,1)$ and let $b$ be the end of $A$ that is approached by both $h_0$ and $h_1$. We let $\sigma(A) = (\sigma(h_0, b, A), \sigma(h_1, b, A))$.

We define some numbers associated with an $H$ and a fixed ordered pair $(h_0, h_1)$ that generates $H$. Let $\alpha(H)$ be the number of orbitals of $H$, and $\beta_i(H)$ be the number of orbitals of $h_i$.

**Remark 2.5.** *We have $\alpha(H) \leq \beta_0(H) + \beta_1(H)$.*

One operation on $H$ will be to replace generators by conjugates of themselves. For an $H$ with a fixed ordered pair of generators $(h_0, h_1)$, we say that a subgroup $G \subseteq H$ is obtained from $H$ by *conjugating generators* if $G$ is generated by the ordered pair $(h'_0, h'_1)$ where $h'_i$ is a conjugate in $H$ of $h_i$ for $i = 1, 2$. We note that if $G = G_k \subseteq G_{k-1} \subseteq \cdots \subseteq G_0 = H$ are given and each $G_{i+1}$ is obtained from $G_i$ by conjugating generators, then $G$ is obtained from $H$ by conjugating generators. We get the following.

**Lemma 2.6.** *Let $H$ be generated by the fixed ordered pair $(h_0, h_1)$. Let the subgroup $G$ of $H$ be obtained from $H$ by conjugating generators. Then*
  *(i) each orbital of $H$ contains at least one orbital of $G$,*
 *(ii) $\alpha(G) \geq \alpha(H)$,*
*(iii) any orbital $A$ of $H$ of pure type $(1,2)$ or $(2,1)$ contains only one orbital $B$ of $G$, and that orbital of $G$ will be of the same pure type as $A$ with $\sigma(B) = \sigma(A)$ when the type is $(2,1)$, and*
*(iv) $\beta_i(G) = \beta_i(H)$ for $i = 1, 2$.*

*Proof.* Item (iv) is immediate. For each orbital $A$ of $H$, some point in $A$ is moved by some $h_i$. Thus any conjugate in $H$ of that $h_i$ moves some point in $A$. This gives (i) which implies (ii). If $A$ is of pure type $(2,1)$ or pure type $(1,2)$, then one of the generators (call it $f$) has $A$ as an orbital. This property is shared by all conjugates in $H$ of $f$ since $A$ is an orbital of $H$. The other generator (call it $g$) approaches one endpoint of $A$ in $A$ and not the other. This property is shared by all conjugates in $H$ of $g$. Lastly, conjugation preserves the sign of an approach, so $\sigma$ is preserved. This proves (iii). $\square$

**Lemma 2.7.** *Let $H$ be generated by the fixed ordered pair $(h_0, h_1)$ and assume that whenever $G \subseteq H$ is obtained from $H$ by conjugating generators, then $\alpha(G) \leq \alpha(H)$. Then $H$ has no orbitals of type $(1,1)$ and all orbitals of types $(1,2)$, $(2,1)$, $(0,2)$ and $(2,0)$ are pure. In addition, if $G \subseteq H$ is obtained from $H$ by conjugating generators, then every orbital of $G$ is an orbital of $H$, and the type, purity, emptiness and $\sigma$ value (if any) with respect to $G$ are the same as with respect to $H$.*

*Proof.* Let an orbital $A$ of $H$ be of type $(2,1)$ and not pure. A conjugate $h_1'$ of $h_1$ has that part of its support in $A$ in a single orbital of $h_0$. If we replace $h_1$ by $h_1'$, then some orbitals of $h_0$ become new pure, empty orbitals of type $(2,0)$ and from (i) of Lemma 2.6, we see that $\alpha(H)$ has strictly increased. Similar arguments take care of orbitals of types $(1,2)$, $(0,2)$ and $(2,0)$. If $H$ has an orbital $A$ of type $(1,1)$ then a conjugate of $h_1$ has that part of its support in $A$ disjoint from the support of $h_0$ and we have another similar argument. The rest of the statement follows easily from Lemma 2.6 and properties of conjugation. $\square$

Another operation will be to replace generators by powers of themselves. Lemma 2.9 will relate to this. One conclusion will be that a wreath product appears, and it turns out that wreath products arise very naturally in our setting. We will need one preliminary lemma about the appearance of wreath products.

We use the term *restricted wreath product* as defined in [13] although we use the symbol $\wr$ to denote the product as opposed to Neumann's "wr." We also use *base*, *top* and *bottom* groups within the structure of a restricted wreath product as defined in [13].

We need one more definition. Given $h$, an orientation preserving self homeomorphism of $I$, and an $x$ not fixed by $h$, we say that a *fundamental domain* of $h$ at $x$ is $x$ together with all points between $x$ and $xh$. Any subset of a fundamental domain of $h$ is disjoint from all of its images under non-zero powers of $h$.

**Lemma 2.8.** *Let $H$ be a group of orientation preserving self homeomorphisms of $I$, and let $h$ be an orientation preserving self homeomorphism of $I$. Assume that the support of $H$ lies in a fundamental domain of $h$. Then $H$ commutes with any conjugate of $H$ by a non-zero power of $h$, and the group $W$ generated by $H$ and $h$ is the restricted wreath product $H \wr \mathbf{Z}$ of $H$ by $\mathbf{Z}$ with $h$ generating the top group $\mathbf{Z}$.*

*Proof.* The first conclusion is immediate. For the second, note that the conjugates of $H$ by powers of $h$ generate the direct sum (not the direct product) of disjoint copies of $H$ indexed over $\mathbf{Z}$, and that this forms a normal subgroup of $W$ that is disjoint from the copy of $\mathbf{Z}$ generated by $h$. The second result now follows from the definition of the restricted wreath product and the characterization of semi-direct products. □

**Lemma 2.9.** *Let $H$ be generated by the fixed ordered pair $(h_0, h_1)$. Let $i$ and $j$ be positive integers and let $G \subseteq H$ be given the fixed ordered pair of generators $(f_0, f_1) = (h_0^i, h_1^j)$. Then the set of orbitals of $H$ and $G$ coincide, and their types, purity, emptiness and $\sigma$ value (if any) with respect to $G$ are the same as with respect to $H$. Further, there is a $K > 0$ so that if $i$ and $j$ are greater than $K$, then:*

 (i) *If $A$ is an orbital of pure type $(2,0)$ or $(0,2)$, then $G_A$ is isomorphic to $\mathbf{Z} \wr \mathbf{Z}$.*
 (ii) *If $A$ is an orbital of pure type $(2,1)$, if $b$ is the endpoint of $A$ approached by both $f_i$ and if $\sigma = \sigma(f_0, b, A)$, then the support of $f_1$ in $A$ is contained in a single orbital of $f_0^\sigma f_1 f_0^{-\sigma}$, and the support of $f_0^{-\sigma} f_1 f_0^\sigma$ in $A$ is contained in a single orbital of $f_1$.*
 (iii) *If $A$ is an orbital of pure type $(2,1)$, if $b$ is the endpoint of $A$ approached by both $f_0$ and $f_1$ in $A$, if $f_1 = f_2$ on a neighborhood of $b$ in $A$, if $\sigma = \sigma(f_1, b, A) = \sigma(f_2, b, A)$, then taking $g'_0 = f_0^\sigma$ to $g_0$ and $g'_1 = f_1^\sigma$ to $g_1$ is an isomorphism from the subgroup of $G_A$ generated by the images in $G_A$ of $(g'_0, g'_1)$ to $F$ as given in Presentation (1).*
 (iv) *If $A$ is an orbital of $G$ of type $(2,2)$, then any fixed point of $f_0$ in $A$ is not a fixed point of either $f_1^{-1} f_0 f_1$ or $f_1 f_0 f_1^{-1}$ and is not a fixed point of any commutator of $f_0^{\pm 1}$ by $f_1^{\pm 1}$, and any fixed point of $f_1$ in $A$ is not a fixed point of either $f_0^{-1} f_1 f_0$ or $f_0 f_1 f_0^{-1}$ and is not a fixed point of any commutator of $f_0^{\pm 1}$ by $f_1^{\pm 1}$.*

*Proof.* That the orbitals of $G$ together with their types and attributes are the same as for $H$ is immediate.

The nature of the claim for $K$ is such that if a value of $K$ is found separately for each of the items on each of the finitely many orbitals, then the maximum value of the different values can be used for the $K$ in the lemma. Thus we will discuss each of the items in the lemma separately for one orbital.

If $A$ is an orbital of $H$ of pure type $(2,0)$, then a sufficiently high power of $h_0$ carries the convex hull of the support of $h_1$ off itself and the rest follows from Lemma 2.8. The argument for pure type $(0,2)$ is symmetric.

If $A$ is an orbital of $H$ of pure type $(2,1)$, then (iii) follows immediately from Lemma 2.1. Item (ii) follows easily because some positive power of $h_0^\sigma$ carries the support of $h_1$ into a single orbital of $h_1$.

We consider (iv). Recall that $A$ is an open set and that we consider only fixed points in $A$. If $x$ is a fixed point of both $f_0$ and a commutator of $f_0^{\pm 1}$ with some $g^{\pm 1}$, then it follows easily that $x$ is also fixed by some conjugate of $f_0^{\pm 1}$ by $g^{\pm 1}$. Thus we only need deal with the conjugates. Since $A$ is an orbital of $H$, we cannot have any points in $A$ fixed by both $h_0$ and $h_1$. All orbits in $A$ under $h_1$ limit (as powers of $h_1$ tend to $\pm\infty$) to fixed points of $h_1$ in $A$ or the endpoints of $A$. Thus no orbit in $A$ under $h_1$ can limit to a fixed point of $h_0$ in $A$. The fixed point set of $f_0$ (which equals the fixed point set of $h_0$) in $A$ consists of a finite number of isolated points and a finite number of closed intervals. Thus there is a $K$ so that $|i| > K$ implies that the image of any fixed point of $f_0$ under $h_1^i$ is disjoint from the fixed point set of $f_0$. This completes half of (iv) and a symmetric argument does the other half. $\square$

*Proof of Theorem 1.* We assume an $H$ as given in the hypothesis of Theorem 1. In the course of the proof, we will replace $H$ by a sequence of groups, each one a subgroup of the preceding. In order to save on notation, we will always refer to the current group as $H$. Since we are trying to find a subgroup of $H$ that is isomorphic to $F$, it will be sufficient to prove that the last version of $H$ is isomorphic to $F$. Because of Remarks 2.2 and 2.3 and Lemma 2.4, we will ignore orbitals on which the projection of $H$ is abelian. In particular, we will ignore pure, empty orbitals of types $(0,2)$ and $(2,0)$. Technically, this is passing to a quotient of $H$, but Lemma 2.4 allows us to safely ignore this distinction.

From Lemma 2.1, we know that there is an orbital $A$ of $H$ and a pair of elements $g_0$ and $g_1$ of $H$ whose restriction to $A$ has the support of $g_1$ in a single orbital $A_0$ of $g_0$, with $A_0$ of pure type $(2,1)$ with respect to the subgroup of $H$ generated by the fixed ordered pair $(g_0, g_1)$. Let $b_0$ be the end of $A_0$ that is approached by both $g_0$ and $g_1$ in $A_0$. By inverting one generator or both if necessary, we can assume that $\sigma(A_0) = (+1, +1)$.

**Step 1.** *Replace $H$ by the subgroup generated by the fixed ordered pair $(g_0, g_1)$.*

We will always use $A_0$ to denote the orbital labeled above as $A_0$ or an orbital with similar properties of a replacement $H$ that is contained in this $A_0$.

At this point, we know that $H$ has only finitely many orbitals.

**Step 2.** *Replace $H$ by a subgroup obtained from $H$ by conjugate replacement so that $\alpha(H)$ is maximized.*

From Remark 2.5 and Lemma 2.6, we know that this can be done and that the resulting group still has an orbital $A_0$ of pure type $(2,1)$ with $\sigma(A_0) = (+1, +1)$. From Lemma 2.7, we know that there are no orbitals of type $(1,1)$ and that all orbitals of types $(0,2)$, $(2,0)$, $(1,2)$ and $(2,1)$ are pure. We use $(g_0, g_1)$ to denote the new fixed ordered pair of generators where each $g_i$ is a conjugate of the previous $g_i$.

**Step 3.** *Replace $H$ by the subgroup generated by powers of the $g_i$ with powers exceeding the value of $K$ as given in Lemma 2.9.*

As before, $(g_0, g_1)$ are the new generators where each $g_i$ is a power of the previous $g_i$. From Lemma 2.9, we know that the orbitals are still as described after the last step, but now we know that the projection of $H$ onto any orbital of type $(2, 0)$ is isomorphic to $\mathbf{Z} \wr \mathbf{Z}$ with $g_0$ generating the top group. Other provisions of Lemma 2.9 will be referred to in the discussion of the next step.

**Step 4.** *Replace $H$ by the subgroup generated by the ordered pair $(g_1, g_0^{-1} g_1 g_0)$.*

Before this step, the projection of $H$ to any orbital of type $(2, 0)$ was isomorphic to $\mathbf{Z} \wr \mathbf{Z}$ with $g_0$ generating the top group, so the projection of the new $H$ to such an orbital of the old $H$ is abelian.

As usual, we recycle notation and use $(g_0, g_1)$ to denote the new pair of generators as given in the statement of the step.

Since orbitals of the previous $H$ of types $(0, 2)$ and $(1, 2)$ were pure, it is elementary that these orbitals are now orbitals of the new $H$ and are of type $(2, 2)$.

If $A$ is an orbital of the previous $H$ of type $(2, 1)$, it is pure, and item (ii) of Lemma 2.9 says that one of the orbitals of the new $H$ in $A$ is also of pure type $(2, 1)$ or pure type $(1, 2)$ and all other orbitals of the new $H$ in $A$ are empty of type $(2, 0)$. We can determine the type and $\sigma$ value of the new non-empty orbital from the following table:

(3)

| old $\sigma$ value | | new type | new $\sigma$ value |
|---|---|---|---|
| $(+1, +1)$ | $\longrightarrow$ | $(2, 1)$ | $(+1, +1)$ |
| $(+1, -1)$ | $\longrightarrow$ | $(2, 1)$ | $(-1, -1)$ |
| $(-1, +1)$ | $\longrightarrow$ | $(1, 2)$ | |
| $(-1, -1)$ | $\longrightarrow$ | $(1, 2)$ | |

In particular, the orbital $A_0$ of the previous $H$ contains an orbital of the new $H$ of pure type $(2, 1)$ that we will continue to call $A_0$ with $\sigma(A_0) = (+1, +1)$.

If $A$ is an orbital of the previous $H$ of type $(2, 2)$, then item (iv) of Lemma 2.9 says that $A$ is also a single orbital of the new $H$ of type $(2, 2)$. Other provisions of that lemma will not be used now.

Note that the only orbitals under consideration at this point are of type $(2, 2)$ and pure types $(2, 1)$ and $(1, 2)$. The orbitals of pure type $(2, 1)$ do not have mixed $\sigma$ values.

**Step 5.** *Repeat Steps 3 and 4 in order two more times.*

From Table (3) and the rest of the discussion in Step 4 we see that after one repetition, there will be orbitals of pure type $(2, 1)$ with $\sigma$ value $(+1, +1)$, pure type $(1, 2)$ and type $(2, 2)$. After the second repetition, there will be orbitals of type $(2, 2)$ and pure type $(2, 1)$ with $\sigma$ value $(+1, +1)$. Among the latter will be the orbital $A_0$.

**Step 6.** *Replace $H$ by the subgroup generated by powers of the $g_i$ with powers exceeding the value of $K$ as given in Lemma 2.9.*

As before, $(g_0, g_1)$ are the new generators where each $g_i$ is a power of the previous $g_i$. From Lemma 2.9, we know that the orbitals are still as described after the last step, so we only have orbitals of type $(2,2)$ and pure type $(2,1)$ with $\sigma$ value $(+1,+1)$.

**Step 7.** *Replace $H$ be the subgroup generated by the ordered pair $(g_1, g_1 g_0^{-1} g_1^{-1} g_0)$.*

We delay renaming the new generators until the end of the discussion of this step.

We start with an orbital $A$ of type $(2,2)$ of the previous $H$. Since the new second generator is a commutator, it is the identity near the end of each orbital of the previous $H$. By item (iv) of Lemma 2.9, the generators of the new $H$ share no fixed points in $A$. Thus any orbital $B$ of the new $H$ in $A$ cannot have both generators approach either end of $B$. Thus all orbitals in $A$ of the new $H$ are of type $(1,1)$, $(2,0)$ or $(0,2)$.

That leaves orbitals of the previous $H$ of pure type $(2,1)$. Let $A$ be such an orbital. Let $b$ be the endpoint of $A$ approached by both $g_0$ and $g_1$. We know that $A$ is an orbital of $g_0$ and that $A$ properly contains the support of $g_1$ that lies in $A$. Item (ii) of Lemma 2.9 says that the pair $(f, h) = (g_1, g_0^{-1} g_1 g_0)$ has an orbital $B$ in $A$ which is of pure type $(2,1)$ with $\sigma$ value $(+1,+1)$ and all other orbitals in $A$ of empty type $(2,0)$. The endpoint $b$ is an endpoint of $B$ and is approached by both $f$ and $h$. The new generators of this step are $(f, fh^{-1})$. We know that $f$ and $h$ have the same slope in $B$ at $b$.

We first study $fh^{-1}$. The only fixed points of $fh^{-1}$ arise from points where $f$ and $h$ agree. Since $f$ has no fixed points in $B$, we cannot get such agreement in $B$ at a fixed point of either $f$ or $h$. Thus the fixed points in $B$ of $fh^{-1}$ are disjoint from any fixed points in $B$ of either $f$ or $h$. Since $f$ and $h$ have the same slope in $B$ at $b$, the support of $fh^{-1}$ is bounded away from $b$ in $B$. Since $B$ is of pure type $(2,1)$, the support of $fh^{-1}$ must approach the end of $B$ opposite that of $b$.

Since the orbital $B$ is of pure type $(2,1)$ for $(f, h)$, the facts in the previous paragraph say that $B$ is also of pure type $(2,1)$ for $(f, fh^{-1})$. Further, we know that $fh^{-1}$ agrees with $f$ on a neighborhood of the end of $B$ opposite $b$. Since $f$ approaches this end negatively, we have $\sigma(B) = (-1, -1)$. All other orbitals of $(f, fh^{-1})$ in $A$ are empty.

We now rename our new generating pair $(g_0, g_1)$. For this pair, we are left with pure orbitals of type $(2,1)$ with $\sigma$ value $(-1,-1)$, one of which is $A_0$, and with orbitals of types $(1,1)$, $(2,0)$ and $(0,2)$.

**Step 8.** *Replace $H$ by a subgroup obtained from $H$ by conjugate replacement so that $\alpha(H)$ is maximized.*

We rename our new generating pair $(g_0, g_1)$.

As in Step 2, we are left with no occurrences of type $(1,1)$ and all types $(2,1)$, $(0,2)$ and $(2,0)$ are pure. Since there was no type $(1,2)$ we end with none. All orbitals of type $(2,1)$ including $A_0$ have $\sigma$ value $(-1,-1)$.

**Step 9.** *Replace each $g_i$ by $g_i^{-1}$.*

This is to convert orbitals of pure type $(2,1)$ to have $\sigma$ value $(+1,+1)$.

**Step 10.** *Replace $H$ by the subgroup generated by equal powers of the $g_i$ with powers exceeding the value of $K$ as given in Lemma 2.9.*

We rename our new generating pair $(g_0, g_1)$.

Ever since Step 7, the orbitals of pure type $(2,1)$ have had the property that the generators have equal slopes at the end of the orbital that they both approach. The point of using equal powers in this step is to preserve that property. Now item (ii) of Lemma 2.9 says that the projection of $H$ to an orbital of type $(2,1)$ is isomorphic to $F$ with the images of the $(g_0, g_1)$ carried to the pair of the same name in Presentation (1). Item (i) of Lemma 2.9 says that the projections to orbitals of types $(2,0)$ and $(0,2)$ are isomorphic to $\mathbf{Z} \wr \mathbf{Z}$.

**Step 11.** *Replace $H$ be the subgroup generated by the ordered pair*
$$(g_0^{-2} g_1 g_0^2 g_1^{-1}, g_0^{-1} g_1 g_0 g_1^{-1}).$$

Since both new generators are commutators, their images in the copies of $\mathbf{Z} \wr \mathbf{Z}$ will commute. We are done when we calculate that the new generators as words in $F$ satisfy the relations of $F$ and do not commute. Since carrying $(g_0, g_1)$ of Presentation (1) to $(g_0, g_1)$ of Presentation (2) is an isomorphism, we can use the relations of Presentation (2) if we define $g_i$ as $g_0^{1-i} g_1 g_0^{i-1}$ for $i > 1$. We let $(G_0, G_1)$ denote the new pair of generators and we define $G_i$ as $G_0^{1-i} G_1 G_0^{i-1}$ for $i > 1$. We must show that $G_1^{-1} G_2 G_1 = G_3$ and $G_1^{-1} G_3 G_1 = G_4$. Our new generators can be rewritten as $(g_3 g_1^{-1}, g_2 g_1^{-1})$. Now

$$\begin{aligned} G_2 &= g_1 g_3^{-1} g_2 g_1^{-1} g_3 g_1^{-1} \\ &= g_1 g_2 g_4^{-1} g_4 g_1^{-1} g_1^{-1} \\ &= g_1 g_2 g_1^{-2} \end{aligned}$$

using $g_3^{-1} g_2 = g_2 g_4^{-1}$ and $g_1^{-1} g_3 = g_4 g_1^{-1}$ which follow directly from the relations in Presentation (2). Similarly, $G_3 = g_1^2 g_2 g_1^{-3}$ and $G_4 = g_1^3 g_2 g_1^{-4}$. Now

$$\begin{aligned} G_1^{-1} G_2 G_1 &= g_1 g_2^{-1} g_1 g_2 g_1^{-2} g_2 g_1^{-1} \\ &= g_1^2 g_2 g_4^{-1} g_4 g_1^{-3} = G_3 \end{aligned}$$

using $g_2^{-1}g_1g_2 = g_1g_2g_4^{-1}$ and $g_1^{-2}g_2 = g_4g_1^{-2}$ which also follow (in two steps) from the relations in Presentation (2). Similarly, one shows that $G_1^{-1}G_3G_1 = G_4$.

To check that $G_0$ and $G_1$ do not commute, we compute a commutator

$$\begin{aligned}G_1G_0^{-1}G_1^{-1}G_0 &= G_1G_2^{-1}\\ &= (g_2g_1^{-1})(g_1^2g_2^{-1}g_1^{-1})\\ &= g_2g_1g_2^{-1}g_1^{-1}\\ &= g_1g_3g_2^{-1}g_1^{-1}\end{aligned}$$

which is trivial if and only if its conjugates are trivial. Conjugating first by $g_1^{-1}$ and then by $g_0^{-1}$ converts this first to $g_3g_2^{-1}$ and then to $g_2g_1^{-1}$ which equals $g_0^{-1}g_1g_0g_1^{-1}$ a commutator of $g_0$ and $g_1$ which cannot be trivial. This completes the proof. □